\documentclass[12pt]{amsart}
 
\usepackage{graphicx} 
\usepackage{amsmath}
\usepackage{amssymb}
\usepackage{mathrsfs}
\usepackage{amscd}
\usepackage{bbm}

\usepackage{enumerate}
\usepackage{multirow}
\usepackage{hyperref}

\newcommand{\calS}{\mathcal S}

\newcommand{\Z}{\mathbb Z}

\newcommand{\be}{\begin{equation}}
\newcommand{\ee}{\end{equation}}

\newcommand{\citep}[1]{\cite{#1}}

\DeclareMathOperator{\lcm}{lcm}

\newtheorem{thm}{Theorem}

\newtheorem{lem}[thm]{Lemma}

\newtheorem{conj}{Conjecture}

\title{Modular periodicity of the Euler numbers and a sequence by Arnold}

\author{Sanjay Ramassamy}
\address{Unit\'e de Math\'ematiques Pures et Appliqu\'ees,
\'Ecole normale sup\'erieure de Lyon,
46 all\'ee d'Italie,
69364 Lyon Cedex 07, France}
\email{sanjay.ramassamy@ens-lyon.fr}

\begin{document}

\maketitle

\begin{abstract}
For any positive integer $q$, the sequence of the Euler up/down numbers reduced modulo $q$ was proved to be ultimately periodic by Knuth and Buckholtz. Based on computer simulations, we state for each value of $q$ precise conjectures for the minimal period and for the position at which the sequence starts being periodic. When $q$ is a power of $2$, a sequence defined by Arnold appears, and we formulate a conjecture for a simple computation of this sequence.
\end{abstract}

\section{Introduction}
\label{sec:introduction}

The sequence of Euler up/down numbers $(E_n)_{n\geq0}$ is the sequence with exponential generating series
\begin{equation}
\sum_{n=0}^{\infty} \frac{E_n}{n!} x^n = \sec x + \tan x.
\end{equation}
It is referenced as sequence A000111 in~\cite{OEIS} and its first terms are
\[
1, 1, 1, 2, 5, 16, 61, 272, 1385, 7936, 50521, 353792, 2702765,\ldots
\]
The numbers $E_n$ were shown by Andr\'e~\cite{A79} to count up/down permutations on $n$ elements (see Section~\ref{sec:powerof2}).

Knuth and Buckholtz~\cite{KB67} proved that for any integer $q\geq1$, the sequence $(E_n \mod q)_{n\geq 0}$ is ultimately periodic. For any $q \geq 1$ we define :
\begin{itemize}
 \item $s(q)$ to be the minimum number of terms one needs to delete from the sequence $(E_n \mod q)_{n \geq 0}$ to make it periodic ;
 \item $d(q)$ to be the smallest period of the sequence $(E_n \mod q)_{n \geq s(q)}$.
\end{itemize}
For example, the sequence $(E_n \mod 3)$ starts with
\[
1,1,1,2,2,1,1,2,2,1,1,2,2,\ldots
\]
so one might expect to have $s(3)=1$ and $d(3)=4$. Clearly $s(1)=0$ and $d(1)=1$. In the remainder of this paper, we formulate precise conjectures for the values of $s(q)$ and $d(q)$ for any $q\geq2$.

\subsection*{Organisation of the paper}
In Section~\ref{sec:oddprime} we reduce the problem to the case when $q$ is a prime power and we conjecture the values of $s(q)$ and $d(q)$ when $q$ is an odd prime power. In Section~\ref{sec:powerof2} we conjecture the values of $s(q)$ and $d(q)$ when $q$ is a power of $2$, after having introduced the Entringer numbers and a sequence defined by Arnold describing the $2$-adic valuation of the Entringer numbers. In Section~\ref{sec:Arnold}, we provide a simple construction which conjecturally yields the Arnold sequence.

\section{Case when $q$ is not a power of $2$}
\label{sec:oddprime}

The following lemma implies that it suffices to know the values of $s(q)$ and $d(q)$ when $q$ is a prime power in order to know the values of $s(q)$ and $d(q)$ for any $q \geq 2$.

\begin{lem}
Fix $q \geq 2$ and write its prime number decomposition as
\begin{equation}
q=\prod_{i=1}^k p_i^{\alpha_i},
\end{equation}
where $k\geq1$, $p_1,\ldots,p_k$ are distinct prime numbers and $\alpha_1,\ldots,\alpha_k$ are positive integers.
Then
\begin{align}
s(q)&=\max_{1 \leq i \leq k} s(p_i^{\alpha_i}) \\
d(q)&=\lcm (d(p_1^{\alpha_1}),\ldots,d(p_k^{\alpha_k})).
\end{align}
\end{lem}

The proof is elementary and uses the Chinese remainder theorem.

When $q$ is an odd prime power, Knuth and Buckholtz~\cite{KB67} found the following :

\begin{thm}[\cite{KB67}]
\label{thm:KB}
Let $p$ be an odd prime number.
\begin{enumerate}
 \item If $p \equiv 1 \mod 4$, then
 \begin{equation*}
 d(p)=p-1.
 \end{equation*} 
 \item If $p \equiv 3 \mod 4$, then
 \begin{equation*}
 d(p)=2p-2.
 \end{equation*}
 \item For any $k\geq1$,
 \begin{equation*}
 s(p^k)\leq k.
 \end{equation*}
 \item For any $k\geq 2$,
 \begin{equation*}
 d(p^k) | p^{k-1} d(p).
 \end{equation*}
\end{enumerate}
\end{thm}

We conjecture the following for the exact values of $s(q)$ and $d(q)$ when $q$ is an odd prime power :

\begin{conj}
\label{conj:oddprime}
Let $p$ be an odd prime number.
\begin{enumerate}
 \item For any $k\geq1$,
 \begin{equation*}
 s(p^k)=k.
 \end{equation*}
 \item For any $k\geq 2$,
 \begin{equation*}
 d(p^k)=p^{k-1} d(p).
 \end{equation*}
\end{enumerate}
\end{conj}
Conjecture~\ref{conj:oddprime} is supported by Mathematica simulations done for all odd prime powers $q<1000$.

\section{Entringer numbers and case when $q$ is a power of $2$}
\label{sec:powerof2}

Formulating a conjecture analogous to Conjecture~\ref{conj:oddprime} for powers of $2$ requires to define, following Arnold~\cite{A91}, a sequence describing the behavior of the $2$-adic valuation of the Entringer numbers.

\subsection{The Seidel-Entringer-Arnold triangle}

\begin{figure}[htpb]
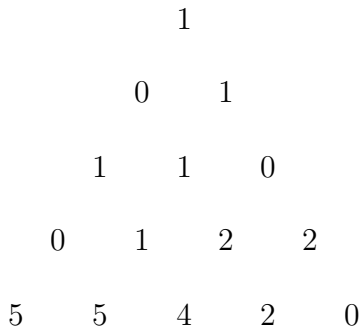

\[
\begin{matrix}
&&&&1&&&& \\
&&&&&&&& \\
&&&0&&1&&& \\
&&&&&&&& \\
&&1&&1&&0&& \\
&&&&&&&& \\
&0&&1&&2&&2& \\
&&&&&&&& \\
5&&5&&4&&2&&0 \\
\end{matrix}
\]
\caption{First five lines of the Seidel-Entringer-Arnold triangle.}
\label{fig:SEAtriangle}
\end{figure}

The Entringer numbers are a refined version of the Euler numbers, enumerating some subsets of up/down permutations. For any $n\geq0$, a permutation $\sigma\in\calS_n$ is called \emph{up/down} if for any $2\leq i \leq n$, we have $\sigma(i-1)<\sigma(i)$ (resp. $\sigma(i-1)>\sigma(i)$) if $i$ is even (resp. $i$ is odd). Andr\'e~\cite{A79} showed that the number of up/down permutations on $n$ elements is $E_n$. For any $1 \leq i \leq n$, the \emph{Entringer number} $e_{n,i}$ is defined to be the number of up/down permutations $\sigma\in\calS_n$ such that $\sigma(n)=i$. The Entringer numbers are usually displayed in a triangular array called the Seidel-Entringer-Arnold triangle, where the numbers $(e_{n,i})_{1 \leq i\leq n}$ appear from left to right on the $n$-th line (see Figure~\ref{fig:SEAtriangle}).

The Entringer numbers can be computed using the following recurrence formula (see for example~\cite{S97}). For any $n\geq2$ and for any $1\leq i \leq n$, we have
\begin{equation}
e_{n,i}=
\begin{cases}
\sum_{j<i} e_{n-1,j} &\text{ if } $n$ \text{ is even} \\
\sum_{j\geq i} e_{n-1,j} &\text{ if } $n$ \text{ is odd}
\end{cases}.
\end{equation}

\subsection{Arnold's sequence}

Replacing each entry of the Seidel-Entringer-Arnold triangle by its $2$-adic valuation, we obtain an infinite triangle denoted by $T$ (see Figure~\ref{fig:2valuationtriangle}).

\begin{figure}[htpb]
\[
\begin{matrix}
&&&&0&&&& \\
&&&&&&&& \\
&&&\infty&&0&&& \\
&&&&&&&& \\
&&0&&0&&\infty&& \\
&&&&&&&& \\
&\infty&&0&&1&&1& \\
&&&&&&&& \\
0&&0&&2&&1&&\infty \\
\end{matrix}
\]
\caption{First five lines of the triangle $T$ of $2$-adic valuations of the Entringer numbers.}
\label{fig:2valuationtriangle}
\end{figure}
We read this triangle $T$ diagonal by diagonal, with diagonals parallel to the left boundary. For any $i\geq1$, denote by $D_i$ the $i$-th diagonal of the triangle $T$ parallel to the left boundary. For example $D_1$ starts with $0,\infty,0,\infty,0,\ldots$. For any $i\geq1$, denote by $m_i$ the minimum entry of diagonal $D_i$. Arnold~\cite{A91} observed that the further away one moves from the left boundary, the higher the $2$-adic valuation of the Entringer numbers becomes. In particular, he observed (without proof) that the sequence $(m_i)_{i\geq1}$ was weakly increasing to infinity. He defined the following sequence : for any $k \geq 1$,
\[
u_k:=\max \left\{i\geq1 | m_i < k \right\}.
\]
In other words, $u_k$ is the number of diagonals containing at least one entry that is not zero modulo $2^k$. The sequence $(u_k)_{k\geq1}$ is referenced as the sequence A108039 in OEIS~\cite{OEIS} and its first few terms are given in Table~\ref{tab:firstterms}.

\begin{table}[htbp]
\centering
\begin{tabular}{|c|c|c|c|c|c|c|c|c|c|c|c|c|c|c|c|c|c|c|c|}
  \hline
  $k$  & 1& 2 & 3 & 4 & 5 & 6 & 7 & 8 & 9 & 10 & 11 & 12 & 13 & 14 & 15 & 16 & 17 & 18 \\
  \hline
  $u_k$ & 2 & 4 & 4 & 4 & 8 & 8 & 8 & 8 & 10 & 12 & 12 & 16 & 16 & 16 & 16 & 16 & 18 & 20 \\
  \hline
\end{tabular}
\caption{The first few values of $u_k$.}
\label{tab:firstterms}
\end{table}

Note that the first few terms given by Arnold were incorrect, because the entry $4$ appeared four times, whereas it should be appearing only three times. We also remark that we cannot define any sequence analogous to $(u_k)$ when studying the $p$-adic valuations of the Entringer numbers for odd primes $p$. Indeed, the $p$-adic valuation $0$ seems to appear in diagonals of arbitrarily high index.

\subsection{Case when $q$ is a power of $2$}

Using the sequence $(u_k)_{k\geq1}$, we formulate the following conjecture for $s(q)$ and $d(q)$ when $q$ is a power of $2$ :

\begin{conj}
\label{conj:powersoftwo}
For any $k\geq 1$, we have
 \begin{equation}
 s(2^k)=u_k.
 \end{equation}
Furthermore, if $k\geq1$ and $k\neq2$, we have
\begin{equation}
d(2^k)= 2^k.
\end{equation}
Finally, we have $d(4)=2$.
\end{conj}
Numerical simulations performed on Mathematica for $k\leq 12$ support Conjecture~\ref{conj:powersoftwo}.

\section{Construction of Arnold's sequence}
\label{sec:Arnold}

In this section we provide a construction which conjecturally yields Arnold's sequence $(u_k)_{k\geq1}$.

We denote by $\Z_+$ the set of nonnegative integers and we denote by
\[
S:=\bigsqcup_{d\geq1}\Z_+^d
\]
the set of all finite sequences of nonnegative integers. We define a map $f:S\rightarrow S$, which maps each $\Z_+^d$ to $\Z_+^{2d}$, as follows. Fix $\underline{x}=(x_1,\ldots,x_d)\in S$. If all the $x_i$'s are equal to $x_d$, we set
\[
f(\underline{x})=(x_d,\ldots,x_d,2x_d,\ldots ,2x_d),
\]
where $x_d$ and $2x_d$ both appear $d$ times on the right-hand side.
Otherwise, define
\[
s:=\max \left\{ 1 \leq i \leq d-1 | x_i \neq x_d \right\}
\]
and set
\[
f(\underline{x})=(x_1,\ldots,x_d,x_1+x_d,\ldots,x_{s-1}+x_d,2x_d,\ldots,2x_d),
\]
where $2x_d$ appears $d-s+1$ times on the right-hand side.
For example, we have
\begin{equation}
f((2,4,4,4))=(2,4,4,4,8,8,8,8)
\end{equation}
and
\begin{equation}
f(2,4,4,4,8,8,8,8)=(2,4,4,4,8,8,8,8,10,12,12,16,16,16,16,16).
\end{equation}
By iterating this function $f$ indefinitely, one produces an infinite sequence :

\begin{lem}
Fix $d\geq1$ and $\underline{x}\in\Z_+^d$. There exists a unique (infinite) sequence $(X_k)_{k\geq 1}$ such that for any $k\geq1$ and for any $n\geq \log_2 (k/d)$, $X_k$ is the $k$-th term of the finite sequence $f^n(\underline{x})$.
\end{lem}

This infinite sequence is called the \emph{$f$-transform} of $\underline{x}$. The lemma follows from the observation that for any $\ell\geq1$ and for any $\underline{y}\in \Z_+^\ell$, $\underline{y}$ and $f(\underline{y})$ have the same first $\ell$ terms.

We can now formulate a conjecture about the construction of the sequence $(u_k)_{k\geq1}$ :
\begin{conj}
\label{conj:Arnold}
Arnold's sequence $(u_k)_{k\geq1}$ is the $f$-transform of the quadruple $(2,4,4,4)$.
\end{conj}
Conjecture~\ref{conj:Arnold} is supported by the estimation on Mathematica of $u_k$ for every $k\leq 512$.

\section*{Acknowledgments}

The author acknowledges the support of the Fondation Simone et Cino Del Duca.

\bibliographystyle{alpha}
\bibliography{bibliographie}

\end{document}